\documentstyle[fleqn,12pt]{article}
\begin{document}
\pagenumbering{arabic}
\setcounter{page}{1}\pagestyle{plain}\baselineskip=20pt

\thispagestyle{empty}
\rightline{MSUMB 96-03, November 1996} \vspace{0.4cm}
\begin{center}
{\bf $h$-deformation of $Gr(2)$}
\end{center}

\vspace{0.5cm}
\noindent
Salih Celik$^*$ \footnote{E-mail: celik@msu.edu.tr}
~and~  Emanullah Hizel$^+$ \footnote{E-mail: hizel@sariyer.cc.itu.edu.tr}

\noindent
$^*$ {\footnotesize Mimar Sinan University, Department of Mathematics, 
80690 Besiktas, Istanbul, TURKEY.}\\
$^+$ {\footnotesize Istanbul Technical University, Department of Mathematics, 
80626 Maslak, Istanbul, TURKEY.}

\vspace{1.5cm}
{\bf Abstract} 

The $h$-deformation of functions on the Grassmann 
matrix group $Gr(2)$ is presented via a contraction of $Gr_q(2)$. 
As an interesting point, we have seen that, in the case of the 
$h$-deformation, both R-matrices of $GL_h(2)$ and $Gr_h(2)$ are 
the same. 

\vfill\eject
In recent years a new class of quantum deformations of Lie groups and 
algebras, the so-called $h$-deformation, has been intensively studied 
by many authors [1-9]. The $h$-deformation of matrix groups can be  
obtained using a contraction procedure. We start with a quantum plane 
and its dual and follow the contraction method of [9]. 

Consider the $q$-deformed algebra of functions on the quantum plane [10] 
generated by $x'$, $y'$ with the commutation rule 
$$x'y' = q y'x'. \eqno(1)$$
Applying a change of basis in the coordinates of the (1) by use of the 
following matrix 
$$g = \left ( \matrix{ 1 & f \cr
                       0 & 1 \cr}
\right) \quad f = {h\over {q - 1}} \eqno(2) $$
one arrives at [9], in the limit $q \rightarrow 1$, 
$$xy = yx + hy^2. \eqno(3)$$
We denote the quantum $h$-plane by $R_h(2)$. 

Similarly, one gets the dual quantum $h$-plane $R^*_h(2)$ as generated 
by $\eta$, $\xi$ with the relations 
$$ \xi^2 = 0 \quad  \eta^2 = h\eta \xi \quad 
   \eta \xi + \xi \eta = 0. \eqno(4) $$

Let 
$$\widehat{A} = \left ( \matrix{ \alpha & \beta \cr
                               \gamma & \delta \cr}\right) $$
be a grassmann matrix in  $Gr(2)$. All matrix elements of 
$\widehat{A}$ are grassmann. We consider linear transformations  
with the following properties:
$$ \widehat{A} : R_h(2) \longrightarrow R_h^*(2) \quad 
   \widehat{A} : R_h^*(2) \longrightarrow R_h(2). \eqno(5)$$
The action on points of $R_h(2)$ and $R^*_h(2)$ of $\widehat{A}$ is 
$$\left(\matrix{ \overline{\eta} \cr \overline{\xi} \cr}\right) = 
  \left(\matrix {\alpha & \beta \cr \gamma & \delta \cr}\right) 
  \left(\matrix{ x \cr y \cr}\right) \quad 
  \left(\matrix{ \overline{x} \cr \overline{y} \cr}\right) = 
  \left(\matrix {\alpha & \beta \cr \gamma & \delta \cr}\right) 
  \left(\matrix{ \eta \cr \xi \cr}\right). \eqno(6)$$
We assume that the entries of $\widehat{A}$ commute 
with the coordinates of $R_h(2)$ and anti-commute with 
the coordinates of $R_h^*(2)$. As a consequence of the linear 
transformations in (5) the vectors 
$\left(\matrix{ \overline{\eta} \cr \overline{\xi} \cr}\right)$ and 
$\left(\matrix{ \overline{x} \cr \overline{y} \cr}\right)$ 
should belong to $R^*_h(2)$ and $R_h(2)$, respectively, which 
impose the following $h$-anti-commutation relations among the 
matrix elements of $\widehat{A}$: 
$$ \alpha \beta + \beta \alpha = h(\alpha \delta + \beta \gamma) \quad 
   \alpha \gamma + \gamma \alpha = 0 $$
$$ \beta  \gamma + \gamma \beta  = h(\delta \gamma - \gamma \alpha) \quad 
   \beta \delta + \delta \beta = -h(\alpha \delta + \gamma \beta) \eqno(7)$$
$$ \alpha \delta + \delta \alpha = h(\gamma \alpha - \delta \gamma) \quad 
   \gamma \delta + \delta \gamma = 0$$
$$ \alpha^2 = -h \gamma \alpha \quad \beta^2 = h(\beta \delta - \alpha \beta + 
   h \alpha \delta) \quad \gamma^2 = 0 \quad \delta^2 = h \delta \gamma.$$
These relations define the $h$-deformation of functions on the grassmann 
matrix group $Gr(2)$, $Gr_h(2)$.

Alternatively, the relations (7) can be obtained by the following similarity  
transformation [9]:
$$\widehat{A}' = g \widehat{A} g^{-1}  \eqno(8)$$
which in our case gives 
$$ \alpha' = \alpha + {h \over {q-1}} \gamma  \quad \beta' = \beta + 
   {h \over {q-1}}( \delta - \alpha - {h \over {q-1}} \gamma)  \eqno(9)$$ 
$$ \hspace*{-5cm} {\gamma' = \gamma \quad \delta' = \delta - {h \over {q-1}} 
    \gamma}.$$ 
and then taking the $q \rightarrow 1$ limit. Here $\alpha'$, $\beta'$, 
$\gamma'$ and $\delta'$ are generators of $Gr_q(2)$, which satisfy the 
following commutation relations [11,12]:
$$ \alpha' \beta' + q^{-1} \beta' \alpha' = 0 \quad 
   \alpha' \gamma' + q^{-1} \gamma' \alpha' = 0 $$
$$ \gamma' \delta' + q^{-1} \delta' \gamma' = 0 \quad 
   \beta' \delta' + q^{-1} \delta' \beta' = 0 \eqno(10)$$
$$ \alpha' \delta' + \delta' \alpha' = 0 \quad 
   \alpha'^2 = \beta'^2 = \gamma'^2 = \delta'^2 = 0 $$
$$ \beta' \gamma' + \gamma' \beta' = (q - q^{-1}) \delta' \alpha'. $$ 
Substituting (9) into (10) one gets the set of relations (7) above. 

The algebra (10) is associative under multiplication and the relations 
in (10) may be also expressed in a tensor product form [11,12] 
$$ {R_q} \widehat{A}_1' \widehat{A}_2' = - 
  \widehat{A}_2' \widehat{A}_1' {R_q} \eqno(11) $$ 
where 
$$ {R_q} = \left ( \matrix{ 
 q + q^{-1} & 0 & 0 & 0 \cr
 0 & 2 & q^{-1} - q & 0 \cr
 0 & q - q^{-1} & 2 & 0 \cr
 0 & 0 & 0 & q + q^{-1} \cr}
\right). \eqno(12) $$ 
Here, since the matrix elements of $\widehat{A}'$ are all grassmann, for 
the conventional tensor products 
$$ \widehat{A}_1' = \widehat{A}' \otimes I ~~\mbox{and}~~ 
\widehat{A}_2' = I \otimes \widehat{A}' \eqno(13) $$
one can write (no-grading)
$$(\widehat{A}_1)^{ij}{}_{{kl}} = 
 \widehat{A}^i{}_{k} \delta^j{}_l \quad 
 (\widehat{A}_2)^{ij}{}_{{kl}} = 
 \delta^i{}_k \widehat{A}^j{}_{l}  \eqno(14)$$
where $\delta$ denotes the Kronecker delta. 
Note that in the limit $q \longrightarrow 1$ the matrix $R_q$ becomes 
twice the 4x4 unit matrix. Notice also that although the algebra (10) is 
an associative algebra of the matrix entries of $\widehat{A}$, $R_q$ 
does not satisfy the quantum Yang-Baxter equation (QYBE) 
$$R_{12} R_{13} R_{23} = R_{23} R_{13} R_{12}. $$
Thus the Yang-Baxter equation is not a necessary 
condition for associativity [see the paragraph after (19) for other 
remarks]. It is obvious that a change of basis in the $R_h(2)$ leads 
to the similarity transformation 
$$ \widehat{A} = g^{-1} \widehat{A}' g  \eqno(15)$$ 
for the quantum grassmann group and the following similarity transformation 
for the corresponding $R$-matrix 
$$ R_{h,q} = (g \otimes g)^{-1} R_q (g \otimes g). \eqno(16)$$ 
If we define the $R$-matrix $R_h$ as 
$$ R_h = \lim_{q \rightarrow 1} R_{h,q}  \eqno(17)$$ 
we get (after dividing by 2) 
$$ {R_h} = \left ( \matrix{ 
 1 & -h & h & h^2 \cr
 0 & 1 & 0 & -h  \cr
 0 & 0 & 1 & h \cr
 0 & 0 & 0 & 1 \cr}
\right). \eqno(18) $$ 
Substituting (9) and (16) into (11) we arrive at the $q \rightarrow 1$ 
limit 
$$ {R_h} \widehat{A}_1 \widehat{A}_2 = - 
  \widehat{A}_2 \widehat{A}_1 {R_h}. \eqno(19) $$ 
An other interesting point is that, although the $R$-matrices of $GL_q(2)$ 
and $Gr_q(2)$ are different, in the case of the $h$-deformation, the 
R-matrices of $GL_h(2)$ and $Gr_h(2)$ are the same [see Ref. 9, 
for the $R$-matrix $R_h$ of $GL_h(2)$]. In the limit 
$q \longrightarrow 1$ both the R-matrices of 
$GL_q(2)$ and $Gr_q(2)$ become the same 4x4 unit matrix. 
Although the R-matrix $R_q$ of $Gr_q(2)$ does not satisfy the QYBE, the 
R-matrix $R_h$ of $Gr_h(2)$ satisfies the QYBE. 

Since the entries of $\widehat{A}$ are all grassmann, a proper inverse 
can not exist. However, the left and right inverses of $\widehat{A}$ 
can be constructed:    
$$\widehat{A}_L^{-1} = \left(\matrix{ 
   \delta + h \gamma & \beta + h \alpha \cr
 - \gamma            & - \alpha \cr}
\right) \eqno(20) $$
$$\widehat{A}_R^{-1} = \left(\matrix{ 
  - \delta & \beta + h \delta \cr
  - \gamma &  \alpha + h \gamma \cr}
\right). \eqno(21) $$

It is now easy to show that 
$$ \widehat{A}^{-1}_{L} \widehat{A} = \Delta_L \eqno(22) $$
$$ \widehat{A} \widehat{A}^{-1}_{R} = \Delta_R \eqno(23)$$
where 
$$\Delta_L = \beta \gamma + \delta \alpha \quad 
  \Delta_R = \gamma \beta + \alpha \delta. \eqno(24)$$
In this case at least formally, $\Delta_L$ and $\Delta_R$ may be 
considered as the left and right quantum (dual) determinants, 
respectively. Note that one can write
$$ \Delta_L \widehat{A}^{-1}_R = \widehat{A}^{-1}_L \Delta_R. \eqno(25)$$

{\it Final Remarks}. We known that all the matrix elements of 
$\widehat{A}$ are grassmann (odd or fermionic) if $\widehat{A}$ is a 
grassmann matrix, i.e., it belongs to $Gr(2)$. Now let $\widehat{A}$ and 
$\widehat{A}'$ be any two anti-commuting (i.e., any element of 
grassmann matrices whose elements $\widehat{A}$ anti-commutes with 
any element of $\widehat{A}'$ ) 
satisfy (10). Then, all the matrix elements of a product 
$A = \widehat{A} \widehat{A}'$ are bosonic (or even) since the elements 
of the matrix product of two grassmann matrices are all bosonic. It can 
also be verified that the matrix elements of $A$ satisfy $q$-commutation 
relations of $GL_q(2)$, i.e., for 
$$A = \left(\matrix{\alpha & \beta \cr \gamma & \delta \cr}\right)
     \left(\matrix{\alpha' & \beta' \cr \gamma' & \delta' \cr}\right) =
\left(\matrix{ a & b \cr c & d \cr}\right)$$
$$ab = q ba \quad ac = q ca \quad bc = cb \eqno(26)$$
{\it etc}. That is, if 
$$\widehat{A}, \widehat{A}' \in Gr_q(2) ~~\Longrightarrow~~ 
  A = \widehat{A} \widehat{A}' \in GL_q(2).$$
In view of these facts, we can say that, there may be no coproduct of the 
form $\Delta(\widehat{A}) = \widehat{A} \dot{\otimes} \widehat{A}$. 
For, this coproduct is invariant under the $q$-commutation relations (26) 
of $GL_q(2)$. These facts also prevent the existence of a coproduct of 
the form $\Delta(\widehat{A}) = \widehat{A}^{t_2} \dot{\otimes} \widehat{A}$ 
where $t_2$ is an involution acting on the elements of $\widehat{A}$. 
Hence a construction of the coproduct along the lines of Ref. 13 is also 
not possible. 

\noindent
{\bf Acknowledgement}

\noindent
This work was supported in part by {\bf T. B. T. A. K.} the 
Turkish Scientific and Technical Research Council. 

\baselineskip=14pt
\def\refname{References}

\end{document}